**Liliane ALFONSI**
*Groupe d'Histoire et de Diffusion des Sciences d'Orsay* (*GHDSO*)
Bât. 407. Université Paris Sud 11
91405 ORSAY *Cédex*
liliane.alfonsi@u-psud.fr

# Un successeur de Bouguer :

# Étienne Bézout (1730–1783) commissaire et expert pour la marine

**Résumé :** Entré à l'Académie Royale des Sciences dès 1758, Étienne Bézout, compte parmi ses tâches, l'obligation d'étudier différents ouvrages et outils présentés à l'Académie. Si dès le début il se voit confier l'étude d'instruments relatifs à la navigation, cela s'amplifiera à partir de 1764, date à laquelle Choiseul, ministre de la Marine, charge Bézout de réorganiser la formation des officiers de ce corps.
Ses visites dans les ports pour faire passer les examens de la Marine, l'amènent à appréhender concrètement les problèmes de la navigation. Cela lui vaudra d'être désigné pour juger la validité des critiques que l'hydrographe Blondeau énonce à l'encontre du traité de Pilotage de Bouguer. Par ailleurs cette situation amène Bézout à être membre de l'Académie de Marine de Brest et à être un des responsables de son association avec l'Académie des Sciences. Nous étudierons ces divers points.
Enfin, Bézout écrit en 1769 un *Traité de navigation* que nous pourrons comparer à l'ouvrage de Bouguer sur le même sujet.

**Mots clés :** Académie des Sciences, Marine, Académie de Marine, Traité de Navigation, Bézout, Bouguer.

**Abstract :** Étienne Bézout, member of the Académie Royale des Sciences, have to study some works and books sended to this institution. In this article, we will look at his responsibility for Navy, before and after 1764, which is the year of Bézout's nomination at the charge of Examinateur des Gardes du Pavillon et de la Marine.
Each year he must go to Brest, Rochefort and Toulon harbours to examine the Gardes de la Marine. This give him titles and qualifications as expert in sailing. Almost in the same time, Étienne Bézout will be member of the Académie de Marine de Brest in 1769. We will see his work in this Academy and as expert, especially his participation to an Academy polemic : Blondeau *versus* Bouguer/Lacaille on a navigation book.
At last, we will study his *Traité de navigation*, written in 1769 and compare it to Bouguer's Navigation book.

Étienne Bézout, né en 1730 à Nemours, fréquente au moins depuis 1755[1] le milieu scientifique parisien et plus particulièrement l'entourage de d'Alembert. Auteur de deux mémoires approuvés par l'Académie royale des sciences - l'un de dynamique en 1756 et l'autre sur les intégrales elliptiques en 1757 -, il est élu à cette même Académie, comme adjoint dans la section de « Méchanique », le 18 mars 1758 (année de la mort de Pierre Bouguer, décédé le 15 août). Comme tous ses collègues académiciens, il a l'obligation d'examiner différents ouvrages et outils présentés à l'Académie. Cette étude a pour but de montrer le rôle de Bézout en tant que commissaire ou expert sur des sujets ayant trait à la navigation, et de façon plus générale, de nous conduire à mieux cerner les implications réciproques entre son statut d'académicien des sciences et son autre statut d'examinateur des écoles de la Marine et d'auteur d'un *Traité de Navigation* pour ces écoles.

---

[1] Liliane Alfonsi, *Étienne Bézout (1730-1783) : mathématicien, académicien et professeur au siècle des Lumières*, Thèse de Doctorat (Paris VI, 2005), 29-34.

### I.  Bézout commissaire et expert pour la Marine

Son rôle de commissaire pour la Marine doit se décomposer en deux parties :
- de 1758 à 1763, années où Bézout est un simple académicien adjoint, section de mécanique ;
- de 1764 à sa mort, période pendant laquelle il fut le responsable des trois écoles des Gardes de la Marine, à Brest, Rochefort et Toulon, pour lesquelles il écrivit un cours de mathématiques dont le dernier volume est un traité de navigation.

**1) de 1758 à 1763**

Le dépouillement systématique des procès verbaux des séances de l'Académie royale des sciences nous donne, pour la période 1758–1763 et pour la marine (certains sujets d'astronomie s'y rapportant), le tableau, placé en annexe, des sujets pour lesquels Bézout fut nommé rapporteur. Nous retiendrons parmi eux : le rapport très négatif sur l'héliomètre de Groumer pour lequel Bézout était le seul commissaire, celui (positif) sur le brise-glaces, qui est de sa main, ce qui, si l'on en croit les habitudes de l'Académie, signifie qu'il est le rapporteur principal ayant vraiment étudié le mémoire.
Bézout eut, durant cette période, 54 mémoires à examiner. Il est intéressant de remarquer que la marine concerne dix d'entre eux, alors que Bézout n'avait encore aucun lien avec elle.

Pendant la séance du mardi 6 septembre 1763, Bézout est élu commissaire pour le prix Rouillé de Meslay [2] que l'Académie veut décerner en 1764 : « L'Académie ayant procédé suivant la forme ordinaire, à la nomination de cinq commissaires pour le prix de 1764, la pluralité des voix a été pour MM. Clairaut, Camus, Lemonnier, Delalande et Bezout. »[3]
Il faut savoir que le choix du sujet et la remise du prix que l'Académie royale des sciences décerne chaque année depuis 1720, sont organisés de la façon suivante : Á l'assemblée publique de Pâques, l'Académie décerne le prix de l'année en cours et les commissaires ayant décernés le prix, proclament le sujet de l'année suivante ou annoncent le report du sujet en cours. Lors de la dernière séance avant les grandes vacances – du 6 septembre au 12

---

[2] Voir sur ce prix, l'article de Guy Boistel dans ce même numéro de la Revue d'Histoire des Sciences.
[3] Archives de l'Académie royale des sciences (Paris) : ARS
Registres (manuscrits) de l'Académie des sciences : *RMAS* (procès-verbaux des séances) 1763, 351.

novembre – on élit les commissaires du prix de l'année suivante, qui reçoivent alors les pièces à examiner. La nomination de Bézout pour le prix de 1764[4], entraîne donc sa participation au choix du sujet pour le prix de 1766[5]. Nous en verrons bientôt les conséquences.

**2) de 1764 à sa mort le 27 septembre 1783**

En 1764 Choiseul, alors ministre de la Marine, décide de réorganiser complètement la formation des officiers de la Marine Royale et charge Étienne Bézout de cette responsabilité. Si la nomination officielle de ce dernier en tant qu'« Examinateur des Écoles des Gardes du Pavillon et de la Marine » est datée du 1$^{er}$ octobre 1764, il est pressenti dès le début de l'année et se rend à Brest au mois de mars pour mieux appréhender sa future fonction. Cela remet en cause son élection comme commissaire du prix, le 21 mars 1764 : « M. Bezout ayant été obligé, par ordre du roi, d'aller à Brest, m'a envoyé par écrit son avis cacheté pour le jugement des pièces du prix dont il est commissaire et sur le sujet à proposer pour 1766. Sur quoi M.Le Président a proposé de délibérer si M.Bezout resterait quoiqu'absent, commissaire du prix ou si on en nommerait un autre à sa place ; il a été décidé qu'il resterait commissaire et que son avis serait remis lors du jugement aux autres commissaires pour être lu et compté comme s'il était présent. »[6]

Bézout reste donc commissaire du prix et son avis va être pris en compte pour le choix du sujet. Or le calcul de la longitude est pleinement d'actualité[7] et parmi les méthodes connues pour ce calcul, celles découlant de l'observation des astres – éclipses de la lune ou des satellites de Jupiter, mesure de la distance d'une étoile à la lune - ont, entre autres sous l'influence de d'Alembert[8], les faveurs de l'Académie. Le sujet est donc à choisir entre l'étude du mouvement de la Lune et celle des satellites de Jupiter[9].

---

[4] Le prix de 1764, qui avait pour sujet « Peut-on expliquer par quelque raison physique pourquoi la lune nous présente toujours à peu près la même face et comment on peut déterminer par les observations et par la théorie, si l'axe de cette planète est sujet à quelque mouvement propre, semblable à celui qu'on connaît dans l'axe de la terre et qui produit la précession des équinoxes et la nutation ? » a été attribué à Lagrange.
[5] Car le sujet de 1763 avait été reporté à 1765.
[6] *RMAS* 1764, 65
[7] Bien que le prix soit pour l'astronomie, il est étroitement lié à la marine par le calcul de la longitude, principal problème de la navigation à cette époque.
[8] D'Alembert pensait que les méthodes d'observations des astres étaient meilleures que les horloges marines. Dans l'article « Longitude » de l'*Encyclopédie* (1765) - où il recommande le traité de navigation de P. Bouguer - il écrit : « Pour découvrir la différence de tems, on s'est servi d'horloges, de montres & d'autres machines, mais toujours en vain, n'y ayant, de tous les instruments propres à marquer le tems, que la seule pendule qui soit assez exacte pour cet effet, & la pendule ne pouvant être d'usage à la mer. D'autres avec des vûes plus saines, & plus de probabilité de succès, vont chercher dans les cieux les moyens de découvrir les *longitudes*. »
[9] Cette dernière est la préférée de d'Alembert, même s'il en reconnaît ses limites en mer. Il écrit, dans le même article « Longitude » : « Cette méthode de déterminer les *longitudes* [satellites de Jupiter] sur terre est aussi exacte qu'on le puisse désirer, & depuis la découverte des satellites de Jupiter, la Géographie a fait de très-grands progrès par cette raison ; [...] Les méthodes qui ont pour fondement des observations de phénomène céleste ayant

À la séance du samedi 7 avril suivant, Lalande intervient pour dire que « les quatre autres commissaires [Lalande, Clairaut, Camus et Le Monnier] s'étoient assemblés et s'étoient trouvés d'accord entre eux et avec M. Bezout pour l'adjudication du prix à une même pièce[10], mais qu'il n'en avoit pas été de même pour la proposition du sujet de 1766. Que deux d'entre eux avoient proposé pour sujet la théorie des satellites de Jupiter et deux autres celle de la Lune. Que l'avis de M. Bezout donnant gain de cause aux premiers, les deux autres avoient refusé de l'admettre par écrit pour la proposition quoiqu'ils l'eussent admis pour l'adjudication »[11], et il demande à l'Académie de statuer sur ce sujet. Après un vote, celle-ci décide, à la majorité, de compter l'avis de Bézout et qu'en conséquence le sujet du prix de 1766 serait la théorie des satellites de Jupiter. Comme l'indique I. Passeron, « la querelle de procédure signale un enjeu scientifique »[12], en effet cet enjeu est celui du calcul de la longitude[13], Clairaut et Lalande étant partisans de la méthode des distances lunaires tandis que d'Alembert préfère celle des satellites de Jupiter. Le vote des commissaires, autres que Bézout, est tenu secret. Mais il est plus que probable que le sujet sur la Lune a été choisi par Clairaut - brouillé avec Le Monnier - et Lalande - ami de Clairaut ainsi qu'ennemi déclaré de d'Alembert et Le Monnier - tandis que les satellites de Jupiter ont été choisis par Bézout et Le Monnier, amis de d'Alembert, ainsi que par Camus, peut-être par amitié pour Bézout.

C'est finalement Lalande qui, malgré son désaccord, rédige le sujet du prix en ces termes : « Quelles sont les inégalités qui doivent s'observer dans les mouvements des quatre satellites de Jupiter à cause de leurs attractions mutuelles, la loi et les périodes de ces inégalités, l'effet qui en résulte sur leurs éclipses et la quantité de ces inégalités selon les meilleures des observations ? Les changements qui paraissent avoir lieu dans les inclinaisons des orbites du second et du troisième satellite doivent surtout être compris dans l'examen de leurs inégalités »[14].

Mais d'Alembert, qui ne fait pas partie des commissaires, se déchaîne contre cette rédaction, dévoilant le rôle qu'il joue dans les coulisses. Il écrit ses *Réflexions sur le programme publié*

---

toutes ce défaut qu'elles ne peuvent être toujours d'usage, parce que les observations ne se peuvent pas faire en tous tems, & étant outre cela d'une pratique difficile en mer, par rapport au mouvement du vaisseau. »

[10] Voir *supra,* note 4.
[11] *RMAS* 1764, 71 r°, v°
[12] Irène Passeron, Une séance à l'Académie au XVIII[e] siècle, *in* Éric Brian et Christiane Demeulenaere-Douyère *Histoire et mémoire de l'Académie des sciences*, Paris : Tec & Doc, 1996, 339-347.
[13] Sur cette question voir : Guy Boistel, *L'astronomie nautique au XVIIIe siècle en France : tables de la lune et longitudes en mer*, thèse de doctorat d'Histoire des sciences et des techniques, Université de Nantes, 2001, Prix André-Jacques Vovard de l'Académie de Marine, 2002, Thèse commercialisée par l'Atelier national de reproduction des Thèses, 2003.
[14] *RMAS* 1764, 98

*au nom de l'Académie pour le prix de 1766*[15], dans lequel les termes du texte du sujet sont presque tous jugés « impropres », voire même « défectueux ». Est alors proposée une addition au texte du sujet : « L'Académie Royale des Sciences de Paris croit devoir avertir les scavants qui travailleront au sujet du prix qu'elle a proposé pour l'année 1766, qu'elle n'entend point exclure l'examen des inégalités que l'action du soleil peut produire dans les mouvements des satellites de Jupiter. » Cette addition signée par Bézout, Camus, Le Monnier et Clairaut, est acceptée par l'Académie[16].

Cette affaire[17] montre l'importance des querelles de personnes et des luttes d'influence à l'Académie, y-compris dans les questions scientifiques.

En dehors de cette polémique Bézout continue son activité habituelle d'académicien et, pendant les années 1764-1765, il est rapporteur pour 19 travaux présentés à l'Académie parmi lesquels six traitent de sujets en rapport avec la marine.

Sur ces six rapports, deux sont de la main de Bézout : sur les horloges marines de Berthoud[18] et sur les satellites de Jupiter de Bailly. Il est aussi nommé commissaire pour la « critique d'un traité de pilotage » de Blondeau, sur lequel nous reviendrons un peu plus loin sur, car le traité critiqué était celui de Pierre Bouguer, édition corrigée par l'abbé de La Caille. Enfin, il faut signaler le cas particulier de la *Connoissance des tems*, publication annuelle de l'Académie des sciences très importante pour la navigation, sur le calendrier et les éphémérides des astres. Le rapport, demandé le 23 novembre 1765 à son sujet, doit être lié à l'opposition entre Lalande, responsable de la revue, et d'Alembert, opposition déjà vue précédemment. En effet, les comptes rendus des séances montrent une polémique visant Lalande : « Sur quelques plaintes qui ont été faites que M. Delalande insérait dans la *Connaissance des temps*, des choses inutiles au but de cet ouvrage, il a été décidé que désormais il serait examiné avant l'impression par MM. De Thury et Bézout, que l'Académie a nommé commissaires à cet effet »[19]. Lalande répond le 1er mars au rapport (introuvable) du 25 janvier 1766, et il obtient le droit d'insérer ses nouvelles tables (parallaxe horizontale et mouvement horaire de la Lune, tables de la réfraction, de la nutation, de l'aberration, etc.), tout en étant obligé de replacer celles qu'il avait ôtées[20].

---

[15] *RMAS* 1764, 192 v°-196 v°
[16] Comme en 1764, c'est Lagrange qui fut le lauréat du prix de 1766.
[17] Voir sur ce sujet Irène Passeron, *op. cit. in* n.12.
[18] Bézout lia au cours du temps des liens d'amitié avec Ferdinand Berthoud, liens concrétisés en 1784 par le mariage de sa propre nièce et héritière avec le neveu de Berthoud, mécanicien de la Marine (voir la thèse de Liliane Alfonsi, *op. cit. in n*. 1, 327).
[19] *RMAS* 1766, 386
[20] Voir sur ce sujet la thèse de Guy Boistel, *op. cit. in* n. 13, 194-205. Les tables des distances lunaires ne seront insérées qu'à partir de 1772. Voir *infra* § IV.

Par ailleurs, depuis qu'il se rend dans les ports pour les examens de la Marine, Bézout en profite pour expérimenter divers instruments que lui confie l'Académie. Ainsi le 17 décembre 1766, cette dernière lui confie une lunette mise au point par l'abbé Alexis-Marie De Rochon pour étudier les satellites de Jupiter. Il la teste à Brest et écrit, le 10 mars 1767, à Duhamel, autre académicien : « Je n'ai pu avoir sur la lunette de M. l'abbé Rochon des expériences aussi complètes que je l'aurai désiré ; et ce n'est ni ma faute, ni la sienne. De 17 jours que j'ai passez à Brest, il n'y en a eu aucun où il n'ait fait un temps affreux. Deux jours seulement, on a pu par intervalles, apercevoir quelques parties de ciel ; c'est l'un de ces deux jours que nous avons fait à terre, l'épreuve de la lunette sur Jupiter. »[21] Bézout pense tout de même que le modèle est valable car il l'a essayé sur une « gabarre » de 4 h 30 à 8 h30 du soir et pour l'observation des satellites, « la lumière de Jupiter est plus que suffisante pour n'être point trop affaiblie par le verre douci »[22]. Il conseille d'aller le tester au Maroc. C'est ce que fera l'abbé Rochon qui partira pour ce pays en avril 1767 sur le vaisseau l'*Union*. Á son retour il présentera à l'Académie des sciences, le 29 août 1767, un mémoire « *Sur des observations faites en mer sur des éclipses par la lune des satellites de Jupiter* » pour lequel Bézout sera nommé commissaire avec d'Alembert, Bailly et Lemonnier. Ils rendent un rapport favorable le 2 septembre suivant, recommandant l'impression dans les Mémoires des Savants étrangers[23].

Sa notoriété scientifique s'impose à la Marine et on lui demande son avis sur des questions d'astronomie rejoignant le problème maritime du calcul de la longitude. C'est ainsi qu'en 1766, Rodier, grand commis à la Marine, le sollicite sur les lettres que le père Esprit Pèzenas (1692-1776) adresse, aussi bien à l'Académie qu'à lui-même, contre Guillaume Saint Jacques de Sylvabelle (1722-1801), ancien élève de Pézenas et nouveau directeur de l'observatoire de Marseille depuis 1763 (à la suite de l'expulsion des jésuites). Pézenas accuse St Jacques de Sylvabelle, de ne pas vouloir lui rendre des effets personnels en prétendant que ces objets appartiennent à l'état[24]. Il l'insulte aussi en tant qu'astronome : « je connois tout le mérite de M. de St Jacques, je sçais qu'il a beaucoup de théorie. Mais je sçais aussi qu'il y a bien loin de la théorie à la pratique. Il a la vüe fort courte et il n'a pas pû trouver dans le ciel aucune trace des deux dernières comètes »[25]. Pézenas fait remarquer que Sylvabelle avait envoyé à l'Académie des sciences, une observation d'éclipse de lune à Marseille, que lui-

---

[21] Archives nationales, à Paris : AN, AN-Marine, G/91, 189
[22] *Ibid*.
[23] Recueils de mémoires écrits par des non-académiciens (d'où « étrangers »), mais publiés par l'Académie royale des sciences
[24] AN-Marine, G/92, 133, lettre du 28 mai 1766
[25] *Ibid*.

même en avait fait autant en observant cette éclipse d'Avignon, ainsi que d'autres astronomes dans d'autres lieux, la sienne et d'autres ont paru dans le Recueil des Savants étrangers de l'Académie, pas celle de Sylvabelle. Il continue sa diatribe : « s'il envoie chaque mois au ministre, toutes les observations qu'il aura faites, sans alléguer le prétexte d'un mauvais temps prétendu, l'Académie pourra aisément décider par la seule inspection et juger de sa capacité. Si elles sont toutes aussi fausses que celle de la dernière éclipse de lune, vous jugerez aisément que votre observatoire est inutile sous la direction d'un tel observateur. Je ne doute point qu'il ne vous envoye dès le premier mois des observations pitoyables, à moins qu'il ne prenne la précaution de les fabriquer dans sa chambre, comme il a fabriqué celle de l'éclipse du soleil en 1764. »[26] Il recommande ensuite de le faire surveiller par M. Mourraille[27], « le seul astronome que je connaisse à Marseille, qui soit en état de faire cette fonction. »[28]

Bézout répond le 7 juin à Rodier qui lui a transmis la lettre de Pézenas, en reprenant, comme à son habitude, méthodiquement et un par un, les arguments de l'abbé pour les contrer successivement :

« Il est bien vrai que Mr de Saint Jacques n'a point déterminé le commencement de l'éclipse du 1$^{er}$ avril 1764, mais M l'abbé Pézenas en sa lettre à M. de l'Isle, pour prouver à celui ci l'impéritie prétendue de Mr. De Saint Jacques, lui fait l'histoire de cette observation : cette histoire prouve que M. Pézenas ignoroit et que M. De Saint Jacques, au contraire, avait connaissance, d'un avertissement que M. de l'Isle avait donné en 1746 aux astronomes au sujet de l'éclipse de cette même année ; c'est pour avoir tenté une expérience à laquelle M. de l'Isle invitait les astronomes, que M. De Saint-Jacques a laissé échapper le commencement de l'éclipse : comme cette expérience était délicate et importante, il n'est pas étonnant que M. De Saint-Jacques n'en ait point retiré le fruit qu'il attendait, et M. de l'Isle n'en est pas surpris non plus : ce qu'il y a de plus surprenant est que M. l'abbé Pézenas ait ignoré la possibilité de ce que M. De Saint-Jacques se proposait et traite tout cela d'absurdité dans sa lettre à M. de l'Isle ; d'ailleurs, Monsieur, n'êtes vous pas surpris de voir l'abbé Pézenas avancer, sans en produire aucune preuve, que M. De Saint-Jacques <u>a fabriqué dans sa chambre les autres phases de cette éclipse</u>. »[29]

---

[26] *Ibid*.
[27] Pierre Mourraille (1721-1808), astronome, mathématicien, élève de Saint-Jacques, secrétaire perpétuel de l'Académie de Marseille de 1768 à 1782 et futur Maire de Marseille entre 1790 et 1793, au caractère violent. Réf : Académie de Marseille, 2003, *Dictionnaire des Marseillais*, Aix-en-Provence, Edisud. Cette précision m'a été communiquée par Guy Boistel et je l'en remercie.
[28] AN-Marine, G/92, 133, lettre du 28 mai 1766
[29] AN-Marine, G/92, 145

Bézout, après avoir défendu Sylvabelle contre les accusations d'incapacité, donne son avis sur Mourraille : « De tous les astronomes de l'Académie, je n'en ai trouvé que deux qui se soient rappelé le nom de M. Mourraille (protégé de M. l'abbé Pézenas), encore a-t-il fallu que je les misse sur la voie pour qu'ils se le rappellassent. C'est sur le rapport du Père Pézenas seulement qu'il est connu de ces deux là, mais il ne l'est de personne par les travaux »[30]

Bézout conteste aussi la destitution de Sylvabelle par l'Académie, réclamée par Pézenas :

« M. l'abbé Pézenas ne pouvant aujourd'hui détruire les éloges qu'il a donnés autrefois à M. De Saint-Jacques, distingue entre la théorie et la pratique et dit qu'il n'a entendu parler que de la première. Mais il serait de la justice de convenir en même tems, qu'en fait d'astronomie, la première est un fort accessoire à la seconde, et qu'un homme dont les ouvrages marquent autant de connoissances et d'adresse dans la théorie, que M. De Saint-Jacques, peut devenir un observateur d'autant plus utile, qu'il a nécessairement tout ce qu'il faut pour bien discuter les observations et les bien employer. Je ne connais point l'état actuel de M. De Saint-Jacques du côté de la pratique ; mais je vois du louche dans les accusations qu'on lui intente, et je serais bien surpris s'il ne parvenait bientôt à faire de bonnes observations. Vous pourrez voir, par la lettre ci jointe, que je vous supplie de vouloir bien faire parvenir à M. l'abbé Pézenas que je ne lui dissimule point ces réflexions. Je suis fâché de voir un homme de son mérite mettre la passion qu'il met dans cette affaire. Je vois l'Académie inondée de ses lettres depuis qu'il s'est mis dans la tête qu'on doit faire juger M. De Saint-Jacques par l'Académie et le déposer sur-le-champ. Il me semble que si l'incapacité de M. De Saint-Jacques était aussi réelle, M. Pézenas devrait penser que l'Académie, si on lui en renvoie le jugement, sçaura bien le connoitre sans être prévenue »[31]

Bézout en arrive alors au prix des objets astronomiques et personnels dont Pézenas demande la restitution ou le remboursement :

« Quant à la désobéissance aux ordres du ministre, dont M. Pézenas accuse M. De Saint-Jacques, je n'ai rien à dire sur ce point. Mais je ne puis m'empêcher de vous observer une contradiction assez singulière entre la lettre que M. Pézenas a eu l'honneur d'écrire à Monseigneur le duc de Praslin, et celle qu'il m'a écrite. Je vois dans la première l'instrument des passages porté à 1200 livres ou <u>aux environs</u>, et dans la dernière à 600 livres seulement. 600 livres sont-ils donc les environs de 1200 livres ? je m'abstiens de commenter cette difficulté. »[32]

---

[30] *Ibid*.
[31] AN-Marine, G/92, 145
[32] *Ibid*.

Puis il conclut dans un *Post-Scriptum* :

« Vous ne serez peut être pas fâché de connoitre comment M. l'abbé Pézenas s'énonçait en 1756 sur le compte de M. De Saint-Jacques qui depuis ce tems est devenu si ignorant. Voici ses paroles : "nous devons le premier (de ces mémoires) à un excellent géomètre, qui connaît Newton comme très peu de personnes le connoissent et à qui on diroit que ce grand homme a transmis son esprit et ses pensées sur le système de l'univers" Préface du 2$^{\text{ème}}$ volume des mémoires rédigés à l'observatoire de Marseille. »[33]

« Par ailleurs j'ai compulsé le plumitif de M. de Fouchy. Il explique que l'Académie n'a pas rejeté son observation [faite par M. De Saint-Jacques] de l'éclipse, mais que, en fait, M. De Saint-Jacques a soumis à l'Académie la solution d'un problème d'astronomie. »[34] Saint-Jacques de Sylvabelle[35] obtiendra gain de cause et ne sera pas destitué.

En 1768, en allant à Brest, Bézout s'arrête à Lorient pour voir des expériences sur le choc des fluides que Thévenard, capitaine dans ce port, veut réaliser. Ce sont alors ses compétences en hydrodynamique qui sont requises par Rodier, lequel lui demande une nouvelle fois son avis. Bézout soutient Thévenard : « Il est constant que pour établir une bonne théorie des fluides, et, par conséquent les principes de la construction des vaisseaux, il faut avoir recours à l'expérience. Il n'est pas moins certain que ces expériences ne peuvent être concluantes qu'autant qu'elles seront faites en grand, qu'elles seront multipliées et exécutées par des personnes intelligentes et exercées. Celles que M. Thévenard médite, me paraissent propres à jeter du jour sur cette partie essentielle à la construction [...] Je pense donc monseigneur, que les fonds que vous voudrez bien lui accorder pour cet objet seront utilement emploïés »[36]. Et Thévenard obtient 5000 livres pour réaliser son projet.

Sans continuer par l'énumération de tous les travaux pour lesquels Bézout a été commissaire pour le ministère de la Marine ou pour l'Académie royale des sciences[37], on notera que sur les 134 rapports auxquels il a participé à l'Académie, 54 ont pour sujet la navigation ou des sujets d'astronomie s'y rattachant pour le calcul de la longitude. Comme il a été rapporteur pour 59 sujets de mathématiques, on peut constater que dans sa tâche d'académicien, les mathématiques et la marine ont eu une importance égale.

---

[33] *Ibid*.]
[34] *Ibid*.
[35] Je remercie Guy Boistel pour les informations suivantes : Cette reconnaissance poussera Saint-Jacques de Sylvabelle à tenter sa chance académique. Le 4 mai 1768, il soumet un *Mémoire sur les longitudes en mer* à l'Académie des sciences (*RMAS*, 1768, 77r°-79v° et pochette séance 11 mai 1768). Examiné par Duséjour et Bailly, ce mémoire est rejeté. Sylvabelle soumet à nouveau son mémoire au jugement de l'Académie de Marine l'année suivante, en décembre 1769. Le mémoire essuie le même échec (Service historique de la Marine à Vincennes (SHMV)-Académie de Marine (AM), t.88, 1). Ce sont les premiers et derniers essais connus de Sylvabelle sur ces questions de navigation (voir Guy Boistel, *op. cit. in n.* 13, Annexe II, 865-866). Il ne produisit plus rien de significatif en astronomie jusqu'à la mise de l'observatoire des Accoules sous tutelle de l'Académie de Marseille en 1782 et restera directeur de l'observatoire jusqu'à sa mort en 1801.
[36] AN-Marine, G/101
[37] On pourra voir pour cela la thèse de Liliane Alfonsi, *op. cit. in n.* 1, 307 et annexe 4.

**II.     La controverse sur le *Nouveau traité de navigation, concernant la théorie et la pratique du pilotage*, de Pierre Bouguer, 2$^e$ édition revue et abrégée par M. l'Abbé de La Caille**

Le vendredi 23 mai 1760, Clairaut et Duhamel du Monceau, académiciens et commissaires nommés par l'Académie royale des sciences pour examiner le *Nouveau traité de navigation, concernant la théorie et la pratique du pilotage*, de Pierre Bouguer, 2$^e$ édition revue et abrégée par M. l'Abbé de La Caille, rendent, en séance, un rapport positif :

« La nouveauté de plusieurs des recherches renfermées dans ces additions & la clarté avec laquelle elles sont toutes exposées, nous on fait penser qu'elles sont très propres à augmenter l'utilité d'un ouvrage dont le public désire une nouvelle édition & qu'elles méritent par conséquent d'être publiées sous le privilège de l'Académie. »[38]

Le traité de Pierre Bouguer, revu par l'Abbé de La Caille est donc édité en 1760. Mais le samedi 11 février 1764, des critiques sur ce traité, émises par Monsieur Blondeau, maître d'hydrographie à Calais, sont soumises au jugement de l'Académie qui nomme commissaires Clairaut et Lemonnier. Ces derniers rendent leur rapport le mercredi 14 mars suivant :

« Nous commissaires nommés par l'Académie, avons examiné des additions et éclaircissements en des cahiers de M. Blondeau, hydrographe du roy à Calais, sur le Traité de Pilotage de M. Bouguer et sur l'édition en abrégé qui en a été faite, il y a quelques années.

Ces sortes d'ouvrages mis entre les mains de ceux qui sçavent à peine les 1$^{ères}$ définitions des élémens de mathématiques, ne sçauraient être trop mis à leur portée et c'est dans les écoles que l'usage que l'on en fait, doit donner assez de vues à ceux qui en font la lecture publique pour éclaircir et détailler plusieurs cas obmis, mais trop utiles dans la pratique. D'ailleurs pour peu qu'un professeur d'hydrographie lise les ouvrages élémentaires avec la sagacité requise, il découvre toujours quelques négligences ou obmissions qu'il est utile de relever et dont les autheurs originaux doivent lui sçavoir gré. C'est ce que M. Blondeau exécute avec soin et qu'il se propose d'imprimer en forme d'additions sans qu'il soit besoin de songer à refondre un ouvrage qui a subi déjà plusieurs degrés de perfection entre les mains des hydrographes qui en ont enrichi le public. Nous avons cru ces additions dignes des éloges de l'Académie et de l'impression. »[39]

---

[38] RMAS 1760, 267-268
[39] ARS, pochette de séance du 14 mars 1764

Ce rapport, qui reste dans la généralité et le flou, approuve les critiques de Blondeau sans que l'on sache de quoi il s'agit, et autorise leur impression. Il ne rejette pas la deuxième édition de l'ouvrage de Bouguer, qu'il recommande de garder inchangée, semblant le considérer comme un ouvrage élémentaire pour débutants.

Cette attitude ne pouvait que choquer le ministre de la Marine, Choiseul, qui, après la destruction de la flotte française pendant la guerre de Sept Ans et la défaite de 1763, s'est attelé à la réorganisation et à la modernisation de la Marine française. Il réagit donc vivement et écrit le samedi 21 juillet 1764 à l'Académie : « Messieurs, J'ai prié M. de Sartine d'arrêter une nouvelle édition du traité de pilotage de feu M. Bouguer publié par feu M. l'abbé de la Caille, sur le rapport qui m'a été fait que cette édition est défectueuse et que M. Blondeau, hydrographe à Calais, y a relevé beaucoup d'erreurs. Je vois par la réponse que me fait M. de Sartine que cette édition a été approuvée par l'Académie au mois de may 1760 et qu'elle a aussi approuvé au mois de mars dernier les observations contradictoires du sieur Blondeau. Comme il importe à la sûreté de la navigation de ne laisser, autant qu'il est possible, aucune incertitude dans ces sortes d'ouvrages, j'ai cru devoir vous prier, Messieurs, d'examiner de nouveau ce traité et les notes de l'hydrographe de Calais afin de décider définitivement. »[40]. Les commissaires nommés pour cette controverse sont Bézout, Camus, Cassini et de Mairan. Ils ne peuvent pas se mettre tout de suite au travail car les critiques de Blondeau semblent avoir été égarées. En effet, le duc de Choiseul écrit deux nouvelles lettres à l'Académie, l'une le 4 août 1764 à Grandjean de Fouchy, secrétaire perpétuel : « J'ay reçu, Monsieur, votre lettre du 23 du mois dernier par laquelle vous m'informez que l'Académie a nommé des commissaires pour examiner la seconde édition du Traité de Navigation de feu M. Bouguer. Il est certain qu'il convient qu'ils examinent en même tems les observations faites sur cet ouvrage par M. Blondeau hydrographe à Calais. J'ay fait chercher ces observations dans mes bureaux et elles ne s'y trouvent pas mais elles pourraient être entre les mains de M. le Comte de Narbonne, inspecteur du dépost des plans et journaux de la Marine auquel l'Académie pourroit les demander ou, s'il ne les avoit pas, elle pourroit faire dire au Sieur Blondeau de luy en envoyer copie. »[41] et comme l'ouvrage de Blondeau n'est pas retrouvé, Choiseul le demande à l'auteur et l'envoie lui-même à l'Académie le 12 novembre, accompagné d'une nouvelle lettre : « J'ai l'honneur de vous envoyer, Messieurs, le mémoire du Sieur Blondeau, hydrographe à Calais, servant de supplément à l'abrégé du Traité de Navigation de M. Bouguer, que j'ay fait demander à cet hydrographe pour être par vous examiné de nouveau

---

[40] *RMAS* 1764, 311
[41] ARS, pochette de séance du 8 août 1764

avec l'édition qu'on veut donner de ce traité. Je vous seroy obligé de me renvoyer ce mémoire lorsque vous en aurez fait usage en me faisant part de votre jugement. »[42]

Cette affaire semble avoir mis très mal à l'aise les quatre commissaires (Bézout, Camus, Cassini, De Mairan) qui en étaient chargés, soit qu'ils n'aient pas été d'accord entre eux ou avec d'autres académiciens - par exemple les premiers nommés sur les remarques de Blondeau : Clairaut et Lemonnier - qu'ils ne voulaient pas affronter, soit qu'ils n'aient pas tous eu d'idée précise sur la valeur des critiques apportées au traité de Bouguer. Ils ne rendent en effet leur rapport que trois ans après, le vendredi 14 août 1767, alors que Choiseul, qui paraissait très attaché à une décision finale dans cette controverse, n'est plus ministre de la Marine (il a laissé la place au Duc de Praslin). De plus ce rapport est introuvable : sa lecture est bien signalée dans le procès verbal de la séance du 14 août, mais il n'a pas été recopié, contrairement aux habitudes de l'Académie, et il ne se trouve pas dans la pochette de séance du jour. En revanche, il est signalé que Lemonnier en a demandé la communication ce jour là, et il est aussi signalé (dans le plumitif de l'année mais pas dans les P.V.), à la séance suivante du 19 août, que ce dernier lit un écrit sur ce sujet, mais la teneur de cet écrit n'est pas révélée. Bézout, Bailly, Duroy et d'Alembert sont nommés rapporteurs pour ces dernières remarques, mais ne rendent jamais leur rapport. Nous pouvons donc émettre l'hypothèse que c'est l'opposition de Lemonnier[43] - précédent commissaire pour Blondeau, nous l'avons dit - qui a fait définitivement enterrer la polémique en entraînant aussi la perte de toute trace des nouvelles conclusions.

On peut aussi penser, avec de sérieuses raisons, que Bézout, lui, était favorable aux critiques de Blondeau, émises, rappelons-le, en mars 1764. En effet, dès sa nomination le 1er octobre 1764 comme Examinateur des Gardes du Pavillon et de la Marine, Bézout fait nommer Blondeau qu'il qualifie de « très bon maître de mathématiques et d'hydrographie à Calais » comme professeur de mathématiques à l'école des Gardes de la Marine de Brest, avec 1500 livres d'appointements[44]. De plus, la grande confiance que Choiseul a toujours montrée envers Bézout (voir thèse de L. Alfonsi, *op. cit. in n.* 1) peut alors expliquer l'acharnement avec lequel ce ministre a tenu au règlement de ce problème, sans toutefois y arriver[45]. Enfin, dans son propre traité de navigation qu'il écrira en 1769, et sur lequel nous reviendrons, Bézout ne suit pas sur certains points le traité de Bouguer, édité par La Caille. Les explications qu'il donne sont peut-être à relier aux critiques de Blondeau.

---

[42] ARS, pochette de séance du 17 novembre 1764
[43] Des informations que Guy Boistel a bien voulu me donner et dont je le remercie, montrent que l'attitude de Lemonnier vis à vis de Blondeau (traiter ses remarques comme des additions mais ne rien vouloir changer au traité Bouguer/Lacaille, comme il l'a écrit avec Clairaut dans le premier rapport) s'explique par le fait que Lemonnier réédite, en 1766, le *Traité de pilotage* de Coubard (voir la thèse de Guy Boistel *op. cit. in n.* 13, 146, 164-165). Il ne peut donc accepter un traité d'un genre nouveau alors qu'il publie un « vieux » traité.
[44] AN-Marine, C/1/3, 67

### III. L'Académie royale de Marine de Brest et son association avec l'Académie royale des sciences de Paris

En 1769, l'Académie de Marine, qui avait vu le jour en 1752 mais qui avait périclité en raison de la guerre de Sept Ans, renaît de ses cendres. Tout naturellement, compte tenu de son rôle dans les écoles de la Marine, on demande à Étienne Bézout d'en faire partie. La nouvelle « Académie Royale de Marine de Brest » se réunit une fois par semaine à partir du 24 mai 1769, date de la 1$^{ère}$ séance, et elle se compose de 10 membres honoraires, 10 membres associés, 20 membres ordinaires et 20 adjoints. On y retrouve le principal fondateur de 1752, Bigot de Morogues, des académiciens des sciences, Bézout, Borda, Bory, Chabert, Duhamel du Monceau, Lalande, Le Monnier, Pingré, Poissonnier, et deux professeurs de mathématiques des Gardes de la Marine de Brest, Blondeau[46] et Duval le Roy. Il y a un seul correspondant, à Orléans[47].

Étienne Bézout qui ne peut participer aux séances que lors de ses passages à Brest, ne sera pas, bien sûr, un des membres les plus actifs. Il envoie tout de même deux écrits à l'Académie, une lettre le 7 octobre 1769[48], relative au défaut de parallélisme des faces des miroirs dans l'usage de l'octant, et un mémoire le 4 février 1772, sur une méthode pour étendre à l'usage des secondes les tables de logarithmes de l'abbé de la Caille[49].

Il est surtout avec Duhamel du Monceau, De Thury, Trudaine et le comte de Maillebois, un des responsables de l'association de l'Académie de Marine avec l'Académie des sciences. L'association est officialisée le 16 mai 1771 avec les règles suivantes :
- Un académicien de la marine ne peut s'intituler des sciences et vice versa ;
- Les officiers de la Marine Royale, membres de l'Académie de Marine, auront droit de séance à l'Académie des sciences quand ils seront à Paris ;
- En dehors de ces officiers, seuls 2 membres de l'Académie de Marine choisis et annoncés, auront ce droit ;

---

[45] Toujours d'après les informations de Guy Boistel, Bézout n'aurait pas osé, à l'Académie, contrer Lemonnier qui était le Préposé officiel au Perfectionnement de la Marine.
[46] Dont la valeur sur les questions maritimes était donc officiellement reconnue.
[47] Voir sur la constitution de l'Académie de Marine, les archives du Service Historique de la Marine de Brest (SHMB), 1/A/108, 237.
[48] SHMV-AM, t. 3, 66, 19 octobre
[49] SHMV-AM, t. 9, 72, 216-220

- Ce droit de séance n'appartiendra pas aux autres membres de l'Académie de Marine, même s'ils ont leur résidence à Paris[50].

L'activité de Bézout à l'Académie de Marine a été très réduite[51] : outre ce qui vient d'être dit, il participera à deux rapports. Le premier, le 22 septembre 1773[52], sur un mémoire d'Eveux de Fleurieu pour un voyage effectué en 1768 et 1769 à bord de l'*Isis*, en vue d'essayer les horloges marines de Berthoud, mémoire à rapprocher d'un autre, présenté celui-ci à l'Académie des sciences le 9 décembre 1767, par Fleurieu aussi, sur les épreuves à faire subir aux horloges de M. Berthoud. Les commissaire étaient Bézout, d'Alembert et Chabert[53].

Le deuxième rapport est écrit le 6 juin 1778, sur le récit d'un voyage effectué par Borda, Pingré et Verdun en 1771 et 1772 à bord de la *Flore*, autour des côtes d'Europe et d'Afrique, pour expérimenter « des méthodes et des instruments mesurant la longitude et la latitude, du bateau et de la côte »[54]. Un rapport a été fait sur le même sujet, le 30 mai 1778 à l'Académie des sciences par Bézout et Bory[55].

On voit sur ces deux exemples le lien entre les deux Académies et le recoupement des rôles de commissaires de Bézout pour ces deux assemblées.

### IV.    Le *Traité de navigation* d'Étienne Bézout

C'est la cinquième et dernière partie (6$^e$ volume) du cours pour la Marine de Bézout. Cet ouvrage, pourtant prévu dès le début de la rédaction de son cours, est publié deux ans après la 4$^e$ partie. Bézout demande des commissaires pour ce *Traité de navigation*, le 16 août 1769, et l'Académie désigne d'Alembert et Duhamel qui rendent leur rapport le 19 août suivant. Ayant reçu l'approbation de l'Académie, son livre est publié avec « Privilège du Roi », dès le 6 décembre 1769. Si Bézout a mis plus de temps pour écrire cette partie, c'est, sans doute, pour recueillir davantage d'informations concrètes et fiables sur la navigation.

En cette année 1769, Bézout, en tant qu'académicien des sciences, a été 19 fois rapporteur sur les sujets les plus divers touchant à la marine. Par ailleurs, ses tournées d'examens dans les ports de Brest, Rochefort et Toulon, lui ont permis d'être informé concrètement des nécessités de la navigation et de réaliser lui-même des expériences, nous l'avons vu. Enfin, sa nomination à l'Académie de Marine de Brest, le 24 avril 1769, soit 4 mois avant la présentation de son traité, lui permet de fréquenter des scientifiques liés à la navigation et des capitaines de vaisseaux reconnus.

---

[50] *RMAS* 1771, 52
[51] Á la décharge d'Étienne Bézout, il faut dire que l'activité de l'Académie de Marine fut très perturbée. Le cahier des présences est là pour montrer l'irrégularité des séances, SHMV-AM, 94, 95, 96, 97.
[52] SHMV-AM, t. 10, 73, 166-177
[53] *RMAS* 1767, 271
[54] SHMV-AM, t. 11, 74, 12-32
[55] *RMAS* 1778, 170

Ce *Traité de Navigation* est composé de quatre sections, les trois premières sont pour tous les élèves, la quatrième ne doit être proposée « qu'à ceux qui veulent se mettre en état de perfectionner l'art de la navigation »

Dans la première section, il définit le cabotage et la navigation hauturière. Cela l'amène à parler de la figure de la Terre, « sphérique ou à très-peu près sphérique ». Bézout définit alors l'équateur, les parallèles, les méridiens, les pôles, les points cardinaux, la latitude et la longitude. Puis il explique les différentes projections existantes pour tracer les cartes. « Après avoir expliqué la construction des cartes marines, il ne s'agit plus, pour être en état d'en faire usage, que de savoir comment on mesure le chemin que fait le navire, & comment on détermine la direction de la route. »[56]

L'usage du « loch »[57], employé pour mesurer la vitesse du navire est expliqué. C'est un sujet que Bézout connaît bien. En effet, en 1764, le ministre de la Marine, Choiseul, lui confie, ainsi qu'à Duhamel du Monceau, en dehors du cadre de l'Académie, l'examen d'un loch, proposé par M. Berly de Blan. C'est Bézout lui-même qui rédige le rapport. Cet écrit[58], montre par ses critiques détaillées – très négatives –, les connaissances de l'auteur.

Vient ensuite l'usage de la boussole pour étudier la « variation », c'est à dire l'angle que fait l'aiguille de la Boussole avec la ligne Nord-Sud. Après beaucoup d'applications pour calculer le chemin du navire, Bézout « fait sentir combien le sillage et la variation, sont incertains ; & combien il est nécessaire d'avoir des moïens de les vérifier et de les rectifier, ce qui donne lieu à la seconde section »[59].

Dans cette deuxième section, « on donne les connoissances d'astronomie utiles aux navigateurs ». Bézout explique comment on détermine la position des astres par rapport à l'écliptique et à l'équateur. Suivent plusieurs exemples traités à l'aide des tables de navigation placées à la fin de l'ouvrage. Bézout termine cette deuxième section en expliquant comment tenir compte des erreurs de parallaxe, et de celles dues à la réfraction de l'air.

La troisième section est celle de l'application des connaissances précédentes grâce aux tables et instruments de navigation. On apprend l'usage du quartier anglais et de l'octant, pour calculer la hauteur d'un astre sur l'horizon. Grâce à ces mesures et aux tables, on obtient la latitude, mais ces déterminations, « quoique bonnes dans la spéculation, ont toutes plusieurs inconvénients dans la pratique surtout à la mer »[60]. Bézout traite plusieurs exemples pour montrer comment le calcul de la latitude permet de corriger la route estimée à l'aide du loch et de la variation. Il donne ensuite les moyens de déterminer l'heure de jour ou de nuit.

---

[56] É. Bézout *Cours de mathématiques à l'usage des gardes du Pavillon et de la Marine*, t. VI, Paris, 1769, p. 26.
[57] Le loch est constitué d'une longue ficelle, divisée en parties égales par des nœuds, et à laquelle est attaché un morceau de bois. On laisse tomber le morceau de bois dans la mer, du côté opposé au vent, la ficelle étant attachée au bateau. Au bout d'une demi-minute, temps mesuré par un sablier, on remonte le morceau de bois en comptant les nœuds de la ficelle. Ce nombre de nœuds donne la vitesse du navire.
[58] Daté du 9 octobre 1764, AN-Marine, G/100, 83-84.
[59] *RMAS* 1769, 314-315

Vient ensuite le problème de la détermination de la longitude, le plus difficile depuis le début de la navigation. Bézout explique les différents moyens possibles, en ne cachant rien de leurs limites :

- Les cartes de la variation de l'aiguille aimantée. Elles donnent les courbes terrestres définies par un même moment magnétique. L'intersection du parallèle déterminé par la latitude avec celle de ces courbes correspondant à l'inclinaison donnée par la boussole, est l'emplacement du navire. Mais, les observations qui ont donné ces courbes ne sont pas assez nombreuses et les courbes varient avec le temps.

- Les montres marines. Elles donnent l'heure du méridien de départ, et la comparaison avec l'heure du lieu du vaisseau, permet d'obtenir la longitude du navire.

Mais l'agitation de la mer et les écarts de température altèrent leurs mouvements. Bézout connaît ce sujet car, à l'Académie des sciences, il avait été nommé[61], le 18 août 1764, commissaire pour l'examen de deux horloges marines[62] l'une de Berthoud et l'autre de Le Roy, les deux spécialistes en France de ces montres. Si l'on ne sait pas qui a rédigé le rapport pour Le Roy, celui de Berthoud est de la main même de Bézout. C'est un rapport favorable, sous réserve d'expérimentation, qui montre le savoir de l'auteur[63].

- Les éclipses de la lune et celles des satellites de Jupiter : Elles « sont visibles au même instant pour tous les lieux où ces astres sont visibles »[64]. Des tables astronomiques donnent les heures de ces éclipses pour un méridien donné. La comparaison de l'heure observée sur le navire avec celle donnée par les tables permet de déterminer la longitude. Mais les tables astronomiques ne sont pas tout à fait exactes et il est très difficile d'observer ces phénomènes en mer de façon satisfaisante. De plus les éclipses de la lune ne peuvent arriver que tous les six mois. Celles des satellites de Jupiter ont lieu presque toutes les nuits, sauf pendant trois mois. Elles seraient donc plus intéressantes pour déterminer la longitude, s'il n'y avait la nécessité d'employer de très longues lunettes. Bézout cite alors le modèle, créé par l'abbé Rochon, qu'il a lui-même testé[65] et dont il pense qu'il est fiable.

---

[60] Étienne Bézout, *op. cit. in n.* 56, 204.
[61] Voir annexe 2.
[62] Il l'avait déjà été le 10 décembre 1763, voir annexe 1.
[63] *RMAS* 1764, 329 et 337.
[64] Étienne Bézout, *op. cit. in n.* 56, 247.

- Une autre méthode est la mesure de la distance d'une étoile à la lune. Elle se fait à l'aide de l'observation par l'octant et des tables donnant les lieux de la lune, disponibles dans le livre de *La Connoissance des Temps* que l'Académie publie chaque année. Cette méthode est celle qui suscite apparemment le moins de critiques.

La quatrième section, en petits caractères, « est destinée à approfondir quelques uns des objets traités dans les trois autres, et à expliquer quelques méthodes qui supposent d'autres connoissances que l'Arithmétique et la Géométrie élémentaire », comme l'écrit d'Alembert dans son rapport[66]. Étienne Bézout explique, entre autres, la correction que l'on doit faire subir aux arcs observés avec l'octant. Il a déjà travaillé sur les problèmes d'optique en les considérant du point de vue géométrique et trigonométrique : en 1767, il a calculé des formules ajoutées au mémoire[67] du Duc de Chaulnes, « Sur quelques expériences relatives à la Dioptrique », mémoire publié dans les *MARS*[68] pour cette même année.

Pour ce cours, il calcule les erreurs que peut occasionner le défaut de parallélisme des miroirs dans l'octant, et il présente, nous l'avons vu, ces résultats à la vingt-deuxième séance de l'Académie de Marine, le 19 octobre 1769, dans une lettre lue par Blondeau[69], dont voici la présentation par ce dernier :

« Tout le monde scait que dans la théorie et dans l'usage de l'octans, on suppose parallèles les faces des miroirs et jusques à présent il paraît qu'on s'est contenté d'exiger ce parallélisme sensiblement. M. Bezout a voulu voir si un petit défaut à cet égard ne produisait pas quelques erreurs de conséquence dans la mesure des arcs. Il a trouvé que, si petit qu'on suppose ce défaut, il peut en résulter une erreur qui n'est point à négliger, il a dressé une table de ces erreurs qu'on trouvera dans le Traité de Navigation qu'il va donner, et ce au moyen d'un excellent quart de cercle. Il a trouvé qu'un très bon octans construit à Paris donne à 85°12' un excès de 5'55''. »[70]

---

[65] Voir *supra* § I-2.
[66] *RMAS* 1769, 314
[67] Voir la thèse de Liliane Alfonsi, *op. cit. in n.* 1, 136-143.
[68] Mémoires de l'Académie Royale des Sciences
[69] Ceci conforte encore l'idée d'une grande proximité de vues entre Bézout et Blondeau sur les problèmes de la Marine et va dans le sens de la thèse exprimée plus haut : Bézout approuvait les critiques de Blondeau à l'égard de la 2ᵉ édition du traité de Bouguer.
[70] SHMV-AM 1769, t. 3, c. 066, 328

Une fois encore[71], Bézout publie dans un cours, une étude et des résultats qu'il aurait pu réserver à une Académie (ici celle de Marine) : la lettre-mémoire envoyée à cette dernière, si elle explique les expériences effectuées, renvoie pour les résultats au *Traité de Navigation*.

Si les grandes lignes du contenu de ce traité sont les mêmes que celles du traité de Bouguer et de la réédition de ce dernier par l'abbé de La Caille, il existe quand même des différences importantes par rapport à ces deux éditions.

La Caille avait déjà traité le calcul de la longitude de façon beaucoup plus conséquente que ne l'avait fait Pierre Bouguer dans son traité de 1753[72]. Mais, bien qu'il passe en revue toutes les méthodes utilisées à l'époque, celle qu'il recommande et qu'il détaille dans son édition est le calcul « par le moyen de la Lune » :

« La méthode de déterminer les longitudes, que nous nous proposons de détailler, se réduit à trois choses :

1°. Á avoir des calculs de la Lune tout faits, selon le modèle que nous donnerons à la fin de ce chapitre. Nous supposons que ces calculs forment un *Almanach nautique*, sur le modèle qui est à la fin de cette Instruction.

2°. Á mesurer sur le vaisseau l'arc de grand cercle, compris dans le ciel entre la Lune & le soleil ou une Étoile choisie.

3°. Á connoître l'instant de temps vrai auquel cet arc de distance aura été mesuré, & de quelle estimation précise on l'aurait trouvé, si la Lune avoit été alors infiniment éloignée de la Terre. Car les calculs tout faits dans l'*Almanach nautique*, indiquent à quel instant on auroit observé sous un Méridien connu, comme celui de Paris, précisément le même arc de distance ; par conséquent la différence de ces deux instants est propre à faire conclure la différence des longitudes entre Paris & le lieu où l'observation aura été faite sur le Navire. »[73] Par la suite, tous ses calculs se réfèrent à cet *Almanach nautique*, qui n'existe pas encore, mais dont il donne – comme il le dit *supra* – un exemple de ce qu'il pourrait être, en donnant, sur une seule page, un modèle de calculs pour juillet 1761[74].

---

[71] Il l'a fait dans son cours d'algèbre (voir la thèse de Liliane Alfonsi, *op. cit. in n*. 1).
[72] Voir la thèse de Guy Boistel *op. cit. in n*. 13.
[73] *Nouveau traité de navigation, contenant la théorie et la pratique du pilotage,* de Pierre Bouguer, 2[e] édition revue et abrégée par M. l'Abbé de La Caille, Paris, 1760, 249.
[74] Inspiré par le modèle de Lacaille, Nevil Maskelyne, astronome anglais, publiera en 1766 le *Nautical Almanach*. En France, ce n'est qu'en 1772 qu'un almanach nautique, comme le concevait Lacaille, sera disponible, publié dans la *Connoissance des tems* par Lalande. Voir la thèse de Guy Boistel *op. cit. in n*. 13, II.1 et III.2, 436 : « Quelle chance fut-ce pour la navigation française si, […] ces astronomes [Lemonnier et Lacaille] avaient uni leurs efforts pour assurer la parution d'un almanach nautique et ce, dès les années 1753-54 ! […] Faisons un rêve ! Un tel almanach aurait très bien pu rassembler les diverses méthodes complémentaires et la réunion des astronomes impliqués dans les progrès de la navigation, Lacaille, Le Monnier, Pingré, Lalande, des géomètres Bouguer, Clairaut, des navigateurs compétents, Mannevillette et Chabert, aurait pu éviter à

C'est sur ce point que réagit Bézout[75]. En parlant de « la méthode de trouver les longitudes par l'observation des distances d'étoiles à la Lune ou au Soleil », il écrit : « Nous nous étions d'abord proposé de suivre, du moins quant au calcul, la méthode que l'on trouve dans l'excellent ouvrage de M. Bouguer (édition de M. l'Abbé de la Caille) ; mais l'Almanach nautique qu'elle suppose, n'existant point ; & n'y ayant pas encore apparence que quelqu'un se charge de sa construction annuelle[76], nous avons cru devoir ne supposer que ce que l'on rencontre plus facilement ; savoir le livre de la *Connoissance des temps*, espèce d'état du ciel, que l'Académie publie chaque année. Mais comme les lieux de la Lune n'y sont calculés que de 12 heures en 12 heures, ce qui n'est pas suffisant pour cet objet, nous avons donné en même temps le moyen d'y suppléer par une règle connue & simple que fournit immédiatement la méthode des interpolations dont nous avons parlé dans l'Algèbre. »[77]

Autre différence dans le traité de navigation de Bézout, sa quatrième section dans laquelle il se propose de « traiter plus à fond plusieurs des objets examinés dans les trois premières »[78]. Il explique comment corriger des erreurs qui se présentent très souvent dans les observations déjà expliquées ou comment améliorer certains calculs :

- corriger l'erreur du défaut de parallélisme des deux miroirs de l'octant (voir *supra*).
- améliorer le calcul de la longitude en mer par la méthode des distances des étoiles à la Lune ou au Soleil[79],

« objets ou nécessaires ou utiles, mais qui n'étant point d'une application indispensable, & exigeant des connaissances ultérieures aux deux premiers volumes de ce cours [*Arithmétique* et *Géométrie*] nous ont paru ne devoir être proposés qu'à ceux qui veulent se mettre en état de perfectionner l'art de la navigation. »[80]

Bien que Bézout ne cite pas Blondeau, on peut se demander, vu leur proximité constatée plusieurs fois précédemment, si ces additions n'ont pas un rapport avec celles qu'avait proposées ce dernier pour améliorer la 2[e] édition du traité de Bouguer par La Caille.

---

l'astronomie nautique française de prendre vingt à trente années de retard sur son homologue britannique. Mais les inerties académiques, le fonctionnement clanique ainsi que le jeu personnalisé de la carrière académique ne favorisaient pas une telle union. Ils ne lui donnaient même aucune chance de se produire. »
[75] Et qu'avait, on peut le penser, réagi Blondeau.
[76] Nous sommes en 1769 donc 9 après l'édition de La Caille. Seul existe le *Nautical Almanach* anglais (voir note 74).
[77] Étienne Bézout, *op. cit. in n.* 56, préface *vj*.
[78] *Ibid*.
[79] Il rend là un hommage à Bouguer puisqu'il écrit en conseillant un nouvel instrument le « Mégamètre » : « Cet instrument, dans la construction duquel M. de Charnières, lieutenant de vaisseau, s'est proposé de rendre l'héliomètre de M. Bouguer applicable à la mesure des distances d'Étoiles à la Lune, a l'avantage de mesurer ces distances avec une précision beaucoup plus grande qu'on ne peut le faire avec l'Octans. » [*Ibid*., 317]
[80] Étienne Bézout, *op. cit. in n.* 56, préface *vij*

## V. Bézout : académicien des sciences et expert pour la marine

Tout au long de cette étude, nous avons vu comment les deux rôles de Bézout, commissaire pour la marine à l'Académie des sciences d'une part, et d'autre part responsable des écoles d'officiers de Marine ainsi que de leur enseignement, ont été étroitement imbriqués avec de fortes implications réciproques.

C'est en tant qu'académicien des sciences qu'il a été amené à connaître les polémiques et les controverses de cette époque sur divers sujets maritimes, à s'intéresser au problème de la longitude en mer, l'une des grandes questions du XVIII$^e$ siècle, et aux perfectionnements apportés à des instruments nécessaires à la navigation. Mais c'est en tant que responsable pour la marine qu'il a pu expérimenter plusieurs de ces instruments, fréquenter les marins et les spécialistes des écoles des Gardes et de l'Académie de Marine de Brest et se forger de cette façon un jugement plus précis sur les besoins pratiques et théoriques de la navigation. Par-là même, il est donc devenu plus performant dans son rôle d'expert pour la marine.

La synthèse de ces deux influences semble être réalisée dans son *Traité de navigation*, dans lequel il utilise à la fois les connaissances qu'il a acquises dans la marine et celle des polémiques découvertes à l'Académie des sciences.

Nous conclurons cette étude sur la façon dont Bézout a été comme Bouguer, mais après lui, commissaire pour la Marine à l'Académie des sciences et auteur d'un traité de navigation, en observant que d'une certaine manière, Bézout a suivi le chemin inverse de Pierre Bouguer. En effet ce dernier, professeur célèbre d'hydrographie au départ, est devenu, à partir de cette fonction, mathématicien et académicien reconnu, alors que Bézout, sans avoir jamais, loin de là, négligé les mathématiques pures, s'est intéressé, les circonstances aidant, à partir de sa fonction de mathématicien et d'académicien aux problèmes de la marine et a contribué, dans l'exercice de ses fonctions, à l'amélioration de leurs solutions.

## ANNEXE

Tableau des sujets en rapport avec la Marine, pour lesquels Bézout fut nommé commissaire pendant la période 1758-1763, la 5$^e$ colonne donnant la date de remise du rapport. Les procès-verbaux ayant été dépouillés jusqu'en 1784, quand cette date n'apparaît pas c'est qu'on ne l'a pas retrouvée jusque là dans ces archives, donc *a priori*, le rapport n'a jamais été rendu.

| Date | Sujet mémoire | Auteur | Autres commissaires | Remise |
|---|---|---|---|---|
| 13/01/1759 | Poussée des voiles | Gobert | Fontaine | |
| 23/08/1760 | Longitudes | abbé Feretti | | |
| 01/07/1761 | Héliomètre | Groumer | | 11/07/1761 |
| 19/05/1762 | Passage de Vénus sur le soleil | Roze | Lemonnier | |
| 20/11/1762 | Barre de Gouvernail | Digard | Clairaut | 27/11/1762 |
| 12/01/1763 | Raideur des Cordages | Thuillier | Clairaut | |
| 16/02/1763 | Brise-glaces | Loriot | Deparcieux, Montigny | 05/03/1763 |
| | Toile à voiles | Daniel | Duhamel, Tillet | 12/03/1763 |
| 14/05/1763 | Toile à voiles (2e mémoire) | Daniel | Duhamel, Tillet | |
| 10/12/1763 | Montre marine | Leroy | Deparcieux, Grandjean, Montigny | |